\documentclass[12pt]{amsart}

\newtheorem{THM}{Theorem}
\newtheorem{LMA}[THM]{Lemma}
\newtheorem{PROP}[THM]{Proposition}
\newtheorem{CORO}[THM]{Corollary}
\newtheorem{CONJ}[THM]{Conjecture}

\numberwithin{equation}{section}

\usepackage{amsmath}
\usepackage{amssymb}
\usepackage{euscript}

\newcommand{\showon}{\begin{eqnarray}}
\newcommand{\showoff}{\end{eqnarray}}

\newcommand{\goesto}{\rightarrow}

\newcommand{\blam}{\boldsymbol{\lambda}}

\newcommand{\A}{\EuScript{A}} 

 \newcommand{\CC}{\mathbb{C}}

\newcommand{\D}{\EuScript{D}} 
\newcommand{\E}{\EuScript{E}} 
\newcommand{\eps}{\varepsilon}

\renewcommand{\H}{\EuScript{H}}

  \newcommand{\NN}{\mathbb{N}}

 \renewcommand{\S}{\EuScript{S}}

\renewcommand{\u}{\mathbf{u}}

\newcommand{\z}{\mathbf{z}}

\begin{document}

\title[Spanning subgraph inequalities]{Enumeration of spanning subgraphs
with degree constraints}
\author{David G. Wagner}
\address{Department of Combinatorics and Optimization\\
University of Waterloo\\
Waterloo, Ontario, Canada\ \ N2L 3G1}
\email{\texttt{dgwagner@math.uwaterloo.ca}}
\thanks{Research supported by the Natural
Sciences and Engineering Research Council of Canada under
operating grant OGP0105392.}

\begin{abstract}
For a finite undirected multigraph $G=(V,E)$ and functions $f,g:V
\goesto\NN$, let $N_{f}^{g}(G,j)$ denote the number of 
$(f,g)$--factors of $G$ with exactly $j$ edges.  The Heilmann-Lieb
Theorem implies that $\sum_{j}N_{0}^{1}(G,j)t^{j}$ is a polynomial
with only real (negative) zeros, and hence that the sequence
$N_{0}^{1}(G,j)$ is strictly logarithmically concave.
Separate generalizations of this theorem were obtained by Ruelle and
by the author.  We unify, simplify, and generalize these
results by means of the Grace--Szeg\"o--Walsh Coincidence Theorem.
\end{abstract}
\maketitle

\section{Introduction.}

By a \emph{graph} $G=(V,E)$ we mean a finite undirected
multigraph.  We identify spanning subgraphs of $G$ with subsets 
$H\subseteq E$ of edges, and let $\deg(H,v)$ denote the degree of 
$v\in V$ in the subgraph $(V,H)$.  The \emph{degree vector} of $H$ is
the function $\deg(H):V\goesto\NN$ given by $\deg(H)(v):=\deg(H,v)$
for all $v\in V$. 
Given functions $f,g:V\goesto\NN$, an \emph{$(f,g)$--factor} is a
spanning subgraph $H$ such that $f\leq \deg(H)\leq g$, in which the
inequalities represent the coordinatewise partial order on $\NN^{V}$.
It is convenient to let natural numbers stand for constant functions on $V$,
so that, for example, a $(0,1)$--factor of $G$ is a matching in $G$.
Let $N_{f}^{g}(G;j)$ denote the number of $(f,g)$--factors of $G$ 
with exactly $j$ edges.  We are concerned here with obtaining 
inequalities among these numbers, and among weighted analogues of them.

This investigation is motivated by the following conjecture.
\begin{CONJ}
For any graph $G=(V,E)$ and functions $f,g:V\goesto\NN$, the sequence
$\{N_f^g(G;j)\}$ is logarithmically concave:\ for all $j$,
$$N_f^g(G;j)^2 \geq N_f^g(G;j-1)N_f^g(G;j+1).$$
\end{CONJ}
The evidence for this is admittedly meagre.  The results reported here --
and the amount of ``slack'' in their derivation -- provide at least some
heuristic support for the conjecture.

Conjecture 1 extrapolates from several results in the literature.
The prototype is the famous (univariate version of a) theorem of
Heilmann and Lieb \cite{HL} (see also Theorem 10.1 of \cite{COSW}).
\begin{THM}[Heilmann--Lieb]
For any graph $G$, the polynomial $\sum_{j}N_{0}^{1}(G;j)t^{j}$ has
only real (strictly negative) zeros.
\end{THM}
Newton's Inequalities (Proposition 12 below) then imply that the 
sequence of coefficients is \emph{strictly logarithmically concave}:
if $N(j)>0$ then $N(j)^{2}>N(j-1)N(j+1)$.

More generally, we have the following (Theorem 3.3 of \cite{Wa1}).
\begin{THM}
For any graph $G=(V,E)$ and functions $f,g:V\goesto\NN$ such that $f\leq 
g\leq f+1$, the polynomial $\sum_{j}N_{f}^{g}(G;j)t^{j}$ has only real 
(nonpositive) zeros.
\end{THM}

Ruelle \cite{Ru} proves a result which relaxes the hypothesis
$f\leq g\leq f+1$.
\begin{THM}[Ruelle]
For a graph $G$ of maximum degree $\Delta\geq 2$, every zero of the
polynomial $\sum_{j}N_{0}^{2}(G;j)t^{j}$ has real part less than
or equal to $-2/\Delta(\Delta-1)^2$.
\end{THM}
A univariate polynomial for which all zeros have negative real part
is said to be \emph{Hurwitz stable}.  This condition implies the
inequalities $N(j)N(j+1)\geq N(j-1)N(j+2)$ among the coefficients
(see Proposition 15 below.)

Another theorem of Ruelle \cite{Ru} involves a weighted version of
the numbers $N_{0}^{2}(G;j)$.  Fix a sequence of nonnegative real
numbers $\u:=\{u_{0},u_{1},u_{2},\ldots\}$ (called \emph{fugacities}),
and for each $H\subseteq E$ let
$$\u_{\deg(H)}:=\prod_{v\in V}u_{\deg(H,v)}.$$
For each natural number $j$ let
$$N(G;\u,j):=\sum_{H\subseteq E:\ \#H=j}\u_{\deg(H)}.$$
For example, when $u_{0}=u_{1}=1$ and $u_{k}=0$ for $k\geq 2$
these are the numbers $N_{0}^{1}(G;j)$.  Similarly, when
$u_{0}=u_{1}=u_{2}=1$ and $u_{k}=0$ for $k\geq 3$ these are the
numbers $N_{0}^{2}(G;j)$.
(More generally, if the functions $f,g$ are not constant on $V$ then the
numbers $N_{f}^{g}(G;j)$ can be expressed in this way only if the fugacities
are \emph{anisotropic}: that is, they are allowed to vary from one vertex to
another.  Our main result is anisotropic but for the purposes of this
introduction the isotropic case above will suffice.)
\begin{THM}[Ruelle]
Let $G$ be a graph with maximum degree $\Delta\geq 1$.  Then for 
fugacities $\u$ satisfying $u_{0}=u_{2}=1$, $u_{1}\geq
\sqrt{2-2/\Delta}$, and $u_{k}=0$ for $k\geq 3$, the polynomial
$\sum_{j}N(G;\u,j)t^{j}$ has only real (strictly negative) zeros. 
\end{THM}

We prove the following result in Section 3.  For a positive
integer $D$ and a sequence of fugacities $\u$, define the generating 
function
$$\Gamma(D,\u,y):=\sum_{k=0}^{D}\binom{D}{k}u_{k}y^{k}.$$
\begin{PROP}
Let $G$ be a graph with maximum degree $\Delta$, and let $D\geq 
\Delta$ be a positive integer.  Let $\u$ be a sequence of fugacities
such that $\Gamma(D,\u,y)$ is a polynomial with only real nonpositive
zeros.  Then the polynomial $\sum_{j}N(G;\u,j)t^{j}$ has only real
nonpositive zeros. 
\end{PROP}
In fact, Proposition 6 can be generalized in two directions -- the fugacities can
be replaced by an anisotropic set $\{\u^{(v)}:\ v\in V\}$ of fugacities,
and the zeros of each $\Gamma(D(v),\u^{(v)},y)$ can be permitted to lie in a
sector centered on the negative real axis (Corollary 19).  The conclusion is
then correspondingly weakened, but often allows the deduction of inequalities
among the coefficients $\{N(G;\{\u^{(v)}\},j)\}$.
Proposition 6 itself implies both Theorems 2 and 5, as is easily verified.
As we shall see, Corollary 19 also implies both Theorem 3 and a wide generalization
of a slight weakening of Theorem 4.  Therefore, Corollary 19 manages to unify
all of the results presented in this introduction.  Moreover, 
its proof is quite straightforward.  Theorem 26 is an application of Corollary 19
which establishes a weak form of Conjecture 1.

Ruelle has a second paper \cite{Ru2} on this subject.  His technique 
uses not only Grace's theorem but also ``Asano contraction'' and some 
intricate geometry of polynomial zeros.  He considers isotropic zero/one
fugacities such that $u_{k}=1$ if and only if  $k\in S$, where $S$ is among
the following the sets: $\{0,1\}$, $\{0,1,2\}$, $\{0,2\}$, $\{0,2,4\}$,
$\{k\ \mathrm{even}\}$, $\{k<\Delta\}$, $\{k\geq 1\}$.  In these cases 
Ruelle produces more precise information about the location of zeros of
$\sum_{j}N(G;\u,j)t^{j}$ than can be obtained by our method.  (On the 
other hand, our method is very easily applied to a wide class of 
vertex degree restrictions.)  Moreover, Ruelle considers factors of
directed graphs in which the indegrees and outdegrees are subject to
separate restrictions.   We indicate briefly at the end how our method
can also be extended to the case of directed graphs.  Perhaps some elaboration
of Ruelle's ideas and our own could lead to further progress towards Conjecture 1.

I thank Alan Sokal for stimulating my interest in the techniques of this paper,
and in particular for showing me the usefulness of the Grace--Szeg\"o--Walsh
Coincidence Theorem.

\section{Preliminaries.}

As indicated in the introduction, we proceed by locating the zeros of
polynomials within certain prescribed sectors and circles.
In this section we first collect the necessary tools, and then explain
the implications for coefficient inequalities.

Let $F(\z)$ be a polynomial in complex variables $\z:=\{z_{v}:\ 
v\in V\}$.  For an open subset $\A\subset\CC$, we say that
$F$ is \emph{$\A$--nonvanishing} if either $F\equiv 0$, or
$z_v\in\A$ for all $v\in V$ implies that
$F(\z)\neq 0$.  If $F\not\equiv 0$ then we say that $F$ is
\emph{strictly} $\A$--nonvanishing.

The following lemma is obvious, and its proof is omitted.
\begin{LMA}
Let $F(\z)$ be $\A$--nonvanishing, and let
$S\subseteq V$. Let $\widetilde{z}_{v}=z_{v}$ if $v\not\in S$, and let
$\widetilde{z}_{v}=z$ if $v\in S$.
Then $F(\widetilde{\z})$ is $\A$--nonvanishing.
\end{LMA}

\begin{LMA}
Let $F_\rho(\z)$ be a family of strictly $\A$--nonvanishing
polynomials indexed by positive real numbers $\rho\in(0,\eps)$.  Assume that
the limit $F(\z):=\lim_{\rho\goesto 0} F_\rho(\z)$
exists.  Then $F$ is $\A$--nonvanishing.
\end{LMA}
\begin{proof}
Each $F_\rho$ is analytic and nonvanishing on the subset $\A^V$ of $\CC^V$.
Since these functions are polynomials, the convergence to $F$ is uniform on
compact subsets of $\CC^{V}$.  By Hurwitz's Theorem, either $F$ is
identically zero or $F$ is nonvanishing on $\A^V$ as well.
\end{proof}

\begin{LMA}
Let $F(\z)$ be $\A$--nonvanishing, and let
$w\in V$.  If $z_{w}$ is fixed at a complex value $\xi_0$ in the closure of $\A$,
then the resulting polynomial in the variables
$\{z_{v}:\ v\in V\smallsetminus\{w\}\}$ is $\A$--nonvanishing.
\end{LMA}
\begin{proof}
Let $\xi:(0,\eps)\goesto\A$ be a continuous function such that $\lim_{\rho\goesto 0}
\xi(\rho)=\xi_0$.  For any positive value of $\rho$, the specialization $z_w=
\xi(\rho)$ results in a polynomial $F_\rho$ that is $\A$--nonvanishing
in the variables $\{z_{v}:\ v\in V\smallsetminus\{w\}\}$.  The result follows from
Lemma 8.
\end{proof}

Of particular importance here are the following subsets of
$\CC$:\\
(i)\ For $0<\theta\leq\pi$, the open sector
$$\S[\theta]:=\{z\in\CC:\ \
z\neq 0\ \mathrm{and}\ |\arg(z)|<\theta\}.$$
(For $z\neq 0$ we use the value of the argument in
the range $-\pi<\arg(z)\leq\pi$.)  For the open right half--plane
$\H:=\S[\pi/2]$ an $\H$--vanishing polynomial is also said to be
\emph{Hurwitz quasi--stable} or have the \emph{half--plane property}.
Notice that if $F$ is $\S[\pi]$--nonvanishing and $F(\z)=0$ then at least one
of the complex numbers $z_{v}$ must be a nonpositive real number.
In particular, a univariate polynomial is $\S[\pi]$--nonvanishing
if and only if it has only real nonpositive zeros.\\
(ii)\ The open unit disc $\D:=\{z\in\CC:\ |z|<1\}$.\\
(iii)\ The open exterior of the unit disc $\E:=\{z\in\CC:\ |z|>1\}$.\\
For a positive real number $\kappa>0$, we let $\kappa\D:=\{z\in\CC:\ |z|<\kappa\}$
and $\kappa\E:=\{z\in\CC:\ |z|>\kappa\}$

The sets $\H$, $\kappa\D$, and $\kappa\E$ are examples of \emph{circular
regions}, in that they are bounded by either circles or straight lines.
For such regions, the Grace--Szeg\"o--Walsh Coincidence Theorem can be
extremely useful.  For a polynomial $F(\z)$, a vertex $w\in V$, and an
integer $D$  greater than or equal to the maximum degree to which $z_{w}$ occurs 
in $F$, the \emph{$D$--th polarization of $z_{w}$ in $F$} is the 
polynomial $P^{D}_{w}F(\z)$ defined as follows. Introduce new
variables $\{z_{w1},\ldots,z_{wD}\}$ and let $e_{k}(z_{w1},\ldots,z_{wD})$
denote the $k$-th elementary symmetric function of 
$\{z_{w1},\ldots,z_{wD}\}$.  Then $P^{D}_{w}F(\z)$ is obtained 
from $F(\z)$ by applying the linear transformation defined by
$z_{w}^{k}\mapsto \binom{D}{k}^{-1}e_{k}(z_{w1},\ldots,z_{wD})$ and linear
extension.
\begin{PROP}[Grace--Szeg\"o--Walsh]
Let $\A$ be a circular region, let $F(\z)$ be $\A$--nonvanishing, let $w\in V$,
and let $D$ be an integer greater than or equal to the maximum degree
$d(w)$ to which $z_{w}$ occurs in $F(\z)$.  If either $D=d(w)$ or
$\A$ is convex then $P^{D}_{w}F(\z)$ is also $\A$--nonvanishing.
\end{PROP}
Theorem 15.4 of Marden \cite{Ma} provides a proof for the case $D=d(w)$.
The theorem also holds for $D>d(w)$ with the additional hypothesis that
$\A$ is convex, as explained in Theorem 2.12 of \cite{COSW}.

\begin{PROP}[Takagi,Weisner]
Let $p(y)=\sum_{k=0}^{n}a_{k}y^{k}$ and $q(y)=\sum_{k=0}^{n}b_{k}y^{k}$
be real polynomials of degree at most $n$.  Assume that $p(y)$ is
is $\S[\pi]$--nonvanishing, and that $q(y)$ is
$\S[\theta]$--nonvanishing for some $\theta>\pi/2$.
Then each of the following polynomials
is $\S[\theta]$--nonvanishing.\\
\textup{(a)}\ $\sum_{k=0}^{n}a_{k}b_{k}y^{k}$;\\
\textup{(b)}\ $\sum_{k=0}^{n}k!a_{k}b_{k}y^{k}$;\\
\textup{(c)}\ $\sum_{k=0}^{n}k!(n-k)!a_{k}b_{k}y^{k}$.
\end{PROP}
For a proof, see Takagi \cite{Ta}, Weisner \cite{We}, or the 
exercises in Section 16 of Marden \cite{Ma}.  Related results can be found in
\cite{GW,W1}.

Location of zeros of generating functions implies combinatorial inequalities
\emph{via} the following results.

\begin{PROP}[Newton's Inequalities]
Let $F(t)=\sum_{j=0}^{d}N(j)t^{j}$ be a univariate polynomial with real
coefficients and only real zeros.  Then
$$\frac{N(j)^{2}}{\binom{d}{j}^{2}}\geq
\frac{N(j-1)}{\binom{d}{j-1}}\cdot\frac{N(j+1)}{\binom{d}{j+1}}$$
for all $1\leq j\leq d-1$.  In particular, $N(j)^{2}>N(j-1)N(j+1)$
whenever $N(j)>0$.
\end{PROP}
See inequality (51) of Hardy--Littlewood--P\'olya \cite{HLP} for a 
proof. 

\begin{PROP}[Aissen--Schoenberg--Whitney]
Let $F(t)=\sum_{j=0}^{d}N(j)t^{j}$ be a univariate polynomial with nonnegative
coefficients.  Then $F(t)$ has only real nonpositive zeros if and only if
every finite minor of the Toeplitz matrix
$$\left[
\begin{array}{ccccc}
N(0) & N(1) & N(2) & N(3) & \cdots\\
0 & N(0) & N(1) & N(2) & \cdots\\
0 & 0 & N(0) & N(1) & \cdots\\
0 & 0 & 0 & N(0) & \cdots\\
\vdots & \vdots & \vdots & \vdots & \ddots
\end{array}
\right]$$
is nonnegative.  In particular, in this case then
$$N(j)^2\geq N(j-1)N(j+1)$$
for all $1\leq j\leq d-1$.
\end{PROP}
See \cite{ASW} or Chapter 8 of Karlin \cite{Ka} for a proof.
 
\begin{PROP}[Folklore]
Let $F(t)=\sum_{j=0}^{d}N(j)t^{j}$ be a univariate polynomial with real
coefficients.  If $F$ is $\S[2\pi/3]$--nonvanishing then
$N(j)^{2}\geq N(j-1)N(j+1)$ for all $1\leq j\leq d-1$.
\end{PROP}
\begin{proof}
Factor $F(t)$ over the reals, and proceed by induction on the degree of
$F(t)$, using the fact that if $p(t)$ and $q(t)$ are polynomials with
logarithmically concave sequences of coefficients, then $p(t)q(t)$ also
has a logarithmically concave sequence of coefficients.
\end{proof}

\begin{PROP}[Hurwitz]
Let $F(t)=\sum_{j=0}^{d}N(j)t^{j}$ be a univariate polynomial with nonnegative
coefficients.  Then $F(t)$ is $\H$--nonvanishing if and only if every minor
of the Hurwitz matrix
$$\left[\begin{array}{ccccccc}
N(1)     & N(3)   & N(5)     & \cdots  &          & 0 \\
N(0)     & N(2)   & N(4)     & \cdots  &          & 0 \\
0        & N(1)   & N(3)     & \cdots  &          & 0 \\
0        & N(0)   & N(2)     & \cdots  &          & 0 \\
\vdots   &        &          & \ddots  &          & \vdots \\
0        &        & \cdots   & N(d-3)  &   N(d-1) & 0 \\
0        & 0      & \cdots   &         &   N(d-2) & N(d)
\end{array}\right]$$
is nonnegative.  In particular, in this case then
$$N(j)N(j+1)\geq N(j-1)N(j+2)$$
for all $1\leq j\leq d-2$., and
$$N(j)^2\geq N(j-2)N(j+2)$$
for all $2\leq j\leq d-2$.
\end{PROP}
See Asner \cite{As} or Kemperman \cite{Ke} for a proof.

\section{Results.}

For a graph $G=(V,E)$, let $\blam:=\{\lambda_{e}:\ e\in E\}$ be positive
real constants indexed by the edges of $G$, and let $\z:=\{z_v:\ v\in V\}$ be
complex variables indexed by the vertices of $G$.  We indicate that the
ends of $e\in E$ are the vertices $v$ and $w$ by writing $vew\in E$
(of course, $v=w$ is possible).  Notice that
$$F(G;\blam,\z):=\prod_{vew\in E}(1+\lambda_{e}z_{v}z_{w})
=\sum_{H\subseteq E}\blam_{H}\z^{\deg(H)}$$
is a weighted multivariate generating function for all spanning subgraphs of $G$,
in which $\blam_{H}:=\prod_{e\in H}\lambda_{e}$.

\begin{PROP}
Let $G=(V,E)$ be a graph, and fix positive constants
$\blam=\{\lambda_e:\ e\in E\}$ such that $\lambda_{\min}\leq\lambda_e\leq
\lambda_{\max}$ for all $e\in E$.  Then\\
\textup{(a)}\  $F(G;\blam,\z)$ is $\H$--nonvanishing.\\
\textup{(b)}\  $F(G;\blam,\z)$ is $\lambda_{\max}^{-1/2}\D$--nonvanishing.\\
\textup{(c)}\  $F(G;\blam,\z)$ is $\lambda_{\min}^{-1/2}\E$--nonvanishing.
\end{PROP}
\begin{proof}
In each case, each factor $1+\lambda_e z_v z_w$ in the product
$F(G;\blam,\z)$ is seen to be nonvanishing in the appropriate region,
from which the result follows.
\end{proof}

The main theorem of this paper is as follows.
Let $G=(V,E)$ be a graph, and for each vertex $v\in V$ fix a sequence
of nonnegative fugacities $\u^{(v)}:=
\{u^{(v)}_{0},u^{(v)}_{1},u^{(v)}_{2},\ldots\}$.
For $H\subseteq E$, let
$$\u_{\deg(H)}:=\prod_{v\in V}u^{(v)}_{\deg(H,v)}.$$
\begin{THM}
Let $G=(V,E)$ be a graph, and let $D:V\goesto\NN$ be such that
$\deg(G)\leq D$. For each $v\in V$, fix
a sequence $\u^{(v)}$ of nonnegative fugacities with generating
function $\Gamma_v(y):=\Gamma(D(v),\u^{(v)},y)$.  Consider the polyomial
$F(G;\{\u^{(v)}\},\z):= \sum_{H\subseteq E}\u_{\deg(H)}\z^{\deg(H)}$.\\
\textup{(a)}\ Fix $0\leq\alpha<\pi/2$. If $\Gamma_v(y)$ is
$\S[\pi-\alpha]$--nonvanishing for all $v\in V$ then $F(G;\{\u^{(v)}\},\z)$ is
$\S[\pi/2-\alpha]$--nonvanishing.\\
\textup{(b)}\ Fix $\kappa>0$. If $\Gamma_v(y)$ is
$\kappa\D$--nonvanishing for all $v\in V$ then $F(G;\{\u^{(v)}\},\z)$ is
$\kappa\D$--nonvanishing.\\
\textup{(c)}\ Fix $\kappa>0$. If $D=\deg(G)$, $u^{(v)}_{D(v)}\neq 0$, and
$\Gamma_v(y)$ is $\kappa\E$--nonvanishing for all $v\in V$ then
$F(G;\{\u^{(v)}\},\z)$ is $\kappa\E$--nonvanishing.
\end{THM}
\begin{proof}
We begin with part (a).
First, we prove the special case in which each $\Gamma(D(v),\u^{(v)},y)$ is
$\S[\pi-\alpha]$--nonvanishing and also has a
nonzero constant term.  Afterward, the general case will be obtained by
a limiting argument.   Thus, in this special case we have, for each $v\in V$,
$$\Gamma_v(y)=\sum_{k=0}^{D(v)}\binom{D(v)}{k}u^{(v)}_{k}
y^{k}=u^{(v)}_0\prod_{i=1}^{D(v)}(1+\xi_{vi}y),$$
for some complex numbers $\{\xi_{vi}\}$ such that either $\xi_{vi}=0$ or
$|\arg(\xi_{vi})|\leq \alpha$ for all $v\in V$ and
$1\leq i\leq D(v)$.

Let $F(G;\z)$ be the polynomial of Proposition 16 in which all 
$\lambda_{e}=1$.  This polynomial is $\H$--nonvanishing.
For each $v\in V$, perform the $D(v)$--fold polarization of the
variable $z_{v}$ in $F(G;\z)$.  Let us denote the resulting
polynomial by  $P^{D}F(G;\{z_{vi}\})$.
Since $D(v)\geq \deg(G,v)$ for each $v\in V$, repeated application of
Proposition 10 implies that $P^{D}F(G;\{z_{vi}\})$ is $\H$--nonvanishing.
Next, replace the variable $z_{vi}$ by $\xi_{vi}z_{vi}$
for all $v\in V$ and $1\leq i\leq D(v)$.  Since either $\xi_{vi}=0$ or
$|\arg(\xi_{vi})|\leq\alpha$ for all $v$ and $i$, and since
$P^{D}F(G;\{z_{vi}\})$ is $\H$--nonvanishing, it follows that
the resulting polynomial
$P^{D}F(G;\{\xi_{vi}z_{vi}\})$ is $\S[\pi/2-\alpha]$--nonvanishing.
Finally, partially diagonalize the variables by making the
substitutions $z_{vi}\mapsto z_{v}$ for all $v\in V$ and
$1\leq i\leq D(v)$.  Repeated application of Lemma 7 shows that
the resulting polynomial $P^{D}F(G;\{\xi_{vi}z_{v}\})$ is still
$\S[\pi/2-\alpha]$--nonvanishing.  Consider the effect of these
operations on $z_{v}^{k}$ in $F(G;\z)$:
$$z_{v}^{k}\mapsto \binom{D(v)}{k}^{-1}e_{k}(\xi_{v1}z_{v},\ldots,
\xi_{vD(v)}z_{v}) = u^{(v)}_{k}z_{v}^{k},$$
since $e_{k}(\xi_{v1},\ldots,\xi_{vD(v)})=\binom{D(v)}{k}u^{(v)}_{k}$.
Therefore,
$$P^{D}F(G;\{\xi_{vi}z_{v}\})
=\sum_{H\subseteq E}\u_{\deg(H)}\z^{\deg(H)}$$
is $\S[\pi/2-\alpha]$--nonvanishing, completing the proof
in this case.

We indicate how the general case of part (a) is derived from the special case 
above by relaxing the hypothesis at one vertex $w\in V$.  Iteration
of this argument for each vertex then yields part (a) as stated.
Accordingly, assume that
$$\Gamma(D(w),\u^{(w)},y)=u^{(w)}_m y^{m}\prod_{i=1}^{D(w)-m}
(1+\xi_{wi}y)$$
for some complex constants $\{\xi_{wi}\}$ with
either $\xi_{wi}=0$ or $|\arg(\xi_{wi})|\leq\alpha$.
Fix a positive real number $\rho>0$, and replace the fugacities
$\u^{(w)}$ by those $\widetilde{\u}^{(w)}$ with generating
function
$$\widetilde{\Gamma}_w(y)=(\rho+ y)^{m}
\prod_{i=1}^{D(v)-m}(1+\xi_{wi}y).$$
One easily verifies that as $\rho\goestoÊ0$,
$$\widetilde{u}^{(w)}_{i}\goesto
\left\{\begin{array}{ll}
0 & \mathrm{if}\ 0\leq i<m,\\
u^{(w)}_{i} & \mathrm{if}\ m\leq i\leq D(w).\end{array}\right.$$
The special case above shows that for any positive value of $\rho>0$, 
$\sum_{H\subseteq E}\widetilde{\u}_{\deg(H)}\z^{\deg(H)}$ is
$\S[\pi/2-\alpha]$--nonvanishing.  By Lemma 8, it follows that
$\sum_{H\subseteq E}\u_{\deg(H)}\z^{\deg(H)}$ is also
$\S[\pi/2-\alpha]$--nonvanishing, as required.

For part (b), the hypothesis implies that every zero of every $\Gamma_v(y)$
has modulus at least $\kappa$.  In particular, for each $v\in V$ we
have $u^{(v)}_0>0$ and
$$\Gamma_v(y)=u^{(v)}_0\prod_{i=1}^{D(v)}(1+\xi_{vi}y)$$
for some complex numbers $\{\xi_{vi}\}$ such that $|\xi_{vi}|\leq 1/\kappa$
for all $1\leq i\leq D(v)$.

Propositions 16 and 10 imply that $P^{D}F(G;\{z_{vi}\})$ is $\D$--nonvanishing.
Replacing each $z_{vi}$ by $\xi_{vi}z_{vi}$, the resulting polynomial
$P^{D}F(G;\{\xi_{vi}z_{vi}\})$ is $\kappa\D$--nonvanishing.
Upon making the substitutions $z_{vi}\mapsto z_{v}$ for all $v\in V$ and
$1\leq i\leq D(v)$, Lemma 7 implies that $F(G;\{\u^{(v)}\},\z)$ is
$\kappa\D$--nonvanishing.

For part (c) the hypothesis implies that every zero of every $\Gamma_v(y)$
has modulus at most $\kappa$, and that $u^{(v)}_{D(v)}>0$.  For each $v\in V$
let $m(v)$ be the multiplicity of $0$ as a root of $\Gamma_v(y)$;\ we may write
$$\Gamma_v(y)=\lim_{\rho\goesto 0}\rho^{m(v)}
u^{(v)}_m(1+y/\rho)^{m(v)}\prod_{i=1}^{D(v)-m(V)}(1+\xi_{vi}y)$$
for some complex numbers $\{\xi_{vi}\}$ such that $|\xi_{vi}|\geq 1/\kappa$
for all $1\leq i\leq D(v)-m(v)$.  Let $\xi_{vi}=1/\rho$ for $D(v)-m(v)+1\leq i\leq D(v)$.

Propositions 16 and 10 imply that $P^{D}F(G;\{z_{vi}\})$ is $\E$--nonvanishing.
For any positive $\rho\leq\kappa$, upon replacing each $z_{vi}$ by $\xi_{vi}z_{vi}$, the
resulting polynomial $P^{D}F(G;\{\xi_{vi}z_{vi}\})$ is $\kappa\E$--nonvanishing.
The substitutions $z_{vi}\mapsto z_{v}$ for all $v\in V$ and
$1\leq i\leq D(v)$, and Lemma 7 imply that $P^{D}F(G;\{\xi_{vi}z_{v}\})$
is $\kappa\D$--nonvanishing.  With $M:=\sum_{v\in V} m(v)$,
the limit of  $\rho^M P^{D}F(G;\{\xi_{vi}z_{v}\})$ as $\rho\goesto 0$ is
$F(G;\{\u^{(v)}\},\z)$.  Lemma 8 then shows that $F(G;\{\u^{(v)}\},\z)$ is
$\kappa\E$--nonvanishing, as desired.
\end{proof}

We can specialize Theorem 17 immediately to obtain a multivariate
generalization of Theorem 3.
\begin{CORO}
For any graph $G=(V,E)$ and functions $f,g:V\goesto\NN$ such that $f\leq 
g\leq f+1$, the polynomial
$$\sum_{H\subseteq E:\ f\leq\deg(H)\leq g}\z^{\deg(H)}$$
is $\H$--nonvanishing.
\end{CORO}
\begin{proof}
Apply Theorem 17(a) by taking $D=\deg(G)$, $\alpha=0$, and
$$\Gamma(D(v),\u^{(v)},y)=\sum_{k=f(v)}^{g(v)}
\binom{D(v)}{k}y^{k}$$
for each $v\in V$.
\end{proof}
The case $f\equiv 0$ and $g\equiv 1$ of Corollary 18 is the multivariate
version of the Heilmann--Lieb theorem \cite{HL}.

With the notation of Theorem 17, for $j\in\NN$ let
$$N(G;\{\u^{(v)}\},j):=\sum_{H\subseteq E:\ \#H=j} \u_{\deg(H)}.$$

\begin{CORO}
Let $G=(V,E)$ be a graph, and let $D:V\goesto\NN$ be such that
$\deg(G)\leq D$.  Fix $0\leq \alpha<\pi/2$. For each $v\in V$, let
$\u^{(v)}$ be a sequence of nonnegative fugacities such that
$\Gamma(D(v),\u^{(v)},y)$ is $S[\pi-\alpha]$--nonvanishing.
Then the univariate polynomial $\sum_{j}N(G;\{\u^{(v)}\},j)t^{j}$ is
$\S[\pi-2\alpha]$--nonvanishing.
\end{CORO}
\begin{proof}
By Theorem 17, the polynomial $\sum_{H\subseteq E}\u_{\deg(H)}\z^{\deg(H)}$
is $\S[\pi/2-\alpha]$--nonvanishing.  Fully diagonalize this, by
making the substitutions $z_{v}\mapsto z$ for all $v\in V$.  By Lemma 7,
the resulting polynomial $\sum_{j}N(G;\{\u^{(v)}\},j)z^{2j}$ 
is also $\S[\pi/2-\alpha]$--nonvanishing.  Since every complex
value $t$ with $|\arg(t)|<\pi-2\alpha$ has a square root $z$ with
$|\arg(z)|<\pi/2-\alpha$, it follows that the polynomial
$\sum_{j}N(G;\{\u^{(v)}\},j)t^{j}$ is
$\S[\pi-2\alpha]$--nonvanishing.
\end{proof}

\begin{CORO}
Let $G=(V,E)$ be a graph, let $D:V\goesto\NN$ be such that
$\deg(G)\leq D$, and to every $v\in V$ assign nonnegative fugacities
$\u^{(v)}$ such that 
$$Q(\u^{(v)},y):=\sum_{k=0}^{D(v)}u^{(v)}_{k}\frac{y^{k}}{k!}$$ is
$\S[\pi-\alpha]$--nonvanishing, with $0\leq\alpha<\pi$.
Then the polynomial
$\sum_{j}N(G;\{\u^{(v)}\},j)t^{j}$ is $(\pi-2\alpha)$--nonvanishing.
\end{CORO}
\begin{proof}
For each $v\in V$, Proposition 11(b) with $p(y)=(1+y)^{D(v)}$ and
$q(y)=Q(\u^{(v)},y)$ shows that the generating function
$\Gamma(D(v),\u^{(v)},y)$ is $\S[\pi-\alpha]$--nonvanishing.
The result follows immediately from Corollary 19.
\end{proof}
Theorem 3.2 of \cite{Wa1} is the case $\alpha=0$ of Corollary 20.
This in turn implies Theorems 2 and 3 above.

Corollary 19 has the following simple consequence.
As grist for the mill we need to locate the zeros of some quadratic
and cubic polynomials within certain sectors.  The proof of Lemma 21 is
a routine calculation, which is omitted.

\begin{LMA}
Let $1\leq k \leq D-1$ be integers, $\beta\geq 0$, let
$$\Gamma(y)=\binom{D}{k-1}y^{k-1}+\beta\binom{D}{k}y^{k}+\binom{D}{k+1}y^{k+1},$$
and let $R=k(D-k)/(k+1)(D-k+1)$.\\
\textup{(a)}  If $\beta\geq \sqrt{2R}$
then $\Gamma(y)$ is strictly $\S[3\pi/4]$--nonvanishing.\\
\textup{(b)}  If $\beta\geq \sqrt{3R}$
then $\Gamma(y)$ is strictly $\S[5\pi/6]$--nonvanishing.\\
\textup{(c)}  If $\beta\geq 2\sqrt{R}$
then $\Gamma(y)$ is strictly $\S[\pi]$--nonvanishing.
\end{LMA}

\begin{LMA}
Let $1\leq k\leq D-2$ be integers, let $\mu\geq 0$, and let
$$\Gamma(y)=\binom{D}{k-1}y^{k-1}+\mu\binom{D}{k}y^{k}+\mu\binom{D}{k+1}y^{k+1}+
\binom{D}{k+2}y^{k+2}.$$
\textup{(a)}  If $\mu\geq 1+\sqrt{2}$
then $\Gamma(y)$ is strictly $\S[3\pi/4]$--nonvanishing.\\
\textup{(b)}  If $\mu\geq 1+\sqrt{3}$
then $\Gamma(y)$ is strictly $\S[5\pi/6]$--nonvanishing.\\
\textup{(c)}  If $\mu\geq 3$
then $\Gamma(y)$ is strictly $\S[\pi]$--nonvanishing.
\end{LMA}
\begin{proof}
Let $p(y)=(1+y)^{D}$ and $q(y)=y^{k-1}(1+\mu y +\mu y^{2} + y^{3})$.
Hypothesis (a) implies that $q(y)$ is $\S[3\pi/4]$--nonvanishing,
hypothesis (b) implies that $q(y)$ is $\S[5\pi/6]$--nonvanishing,
and hypothesis (c) implies that $q(y)$ is $\S[\pi]$--nonvanishing.
The conclusion now follows from Proposition 11(a).
\end{proof}
Lemma 22 ignores some ``finite $D$'' effects that can be significant,
especially for small $D$.  The general case is quite complicated, but
Lemma 23 is indicative of the possible improvement.  Again, the
elementary calculation is omitted.
\begin{LMA}
Let $p\geq 1$ be an integer, $D=2p+1$, $\mu\geq 0$, and let
$$\Gamma(y)=
\binom{2p+1}{p-1}y^{p-1}+\mu\binom{2p+1}{p}y^{p}+\mu\binom{2p+1}{p+1}y^{p+1}+
\binom{2p+1}{p+2}y^{p+2}.$$
\textup{(a)}  If $\mu\geq (1+\sqrt{2})p/(p+2)$ then
$\Gamma(y)$ is strictly $\S[3\pi/4]$--nonvanishing.\\
\textup{(b)}  If $\mu\geq (1+\sqrt{3})p/(p+2)$ then
$\Gamma(y)$ is strictly $\S[5\pi/6]$--nonvanishing.\\
\textup{(c)}  If $\mu\geq 3p/(p+2)$ then
$\Gamma(y)$ is strictly $\S[\pi]$--nonvanishing.
\end{LMA}

\begin{PROP}
Let $G=(V,E)$ be a graph, let $D:V\goesto\NN$ be such that
$\deg(G)\leq D$, and assign nonnegative fugacities $\u^{(v)}$
to each $v\in V$ so that each $\Gamma(D(v),\u^{(v)},y)$ has
$r(v)\leq 4$ nonzero terms, of consecutive degrees.  In each case
below, assume that the generating functions for $v\in V$ with
$r(v)\geq 3$ satisfy the given hypotheses.\\
\textup{(a)} Lemmas 21, 22, and 23, part \textup{(a)}:\ then
$F(G;\{\u^{(v)}\},t)$ is $\H$--nonvanishing.  In this case
the inequalities of Proposition 15 hold for $\{N(G;\{\u^{(v)}\},j)\}$.\\
\textup{(b)} Lemmas 21, 22, and 23, part \textup{(b)}:\ then
$F(G;\{\u^{(v)}\},t)$ is $\S[2\pi/3]$--nonvanishing. In this case
the inequalities of Proposition 14 hold for $\{N(G;\{\u^{(v)}\},j)\}$.\\
\textup{(c)}  Lemmas 21, 22, and 23, part \textup{(c)}:\ then
$F(G;\{\u^{(v)}\},t)$ has only real nonpositive zeros. In this case
the inequalities of Propositions 12 and 13 hold for $\{N(G;\{\u^{(v)}\},j)\}$.
\end{PROP}
\begin{proof}
In part (a) every $\Gamma(D(v),\u^{(v)},y)$ is $\S[3\pi/4]$--nonvanishing.
In part (b) every $\Gamma(D(v),\u^{(v)},y)$ is $\S[5\pi/6]$--nonvanishing.
In part (c) every $\Gamma(D(v),\u^{(v)},y)$ is $\S[\pi]$--nonvanishing.
The result follows immediately from Corollary 19.
\end{proof}

Of most interest combinatorially is the case in which all fugacities are
either zero or one.  All we obtain in this direction is the following rather
limited result.
\begin{PROP}
Let $G=(V,E)$ be a graph, and let $f,g:V\goesto\NN$ be functions
such that $f\leq g\leq f+2$ and $g\leq\deg(G)$.\\
\textup{(a)}  Then the polynomial $\sum_{j}N_{f}^{g}(G;j)t^{j}$
is $\S[\pi/3]$--nonvanishing.\\
\textup{(b)}  Assume furthermore that for every $v\in V$, either
$g(v)\leq f(v)+1$, or $f(v)=0$, or $g(v)=\deg(G,v)$, or 
$\deg(G,v)\leq 5$.  Then the polynomial
$\sum_{j}N_{f}^{g}(G;j)t^{j}$ is $\H$--nonvanishing.
Thus the inequalities of Proposition 15 hold for $\{N_f^g(G;j)\}$.
\end{PROP}
\begin{proof}
We apply Corollary 19 with $D=\deg(G)$ and fugacities given by
$$u_{i}^{(v)}:=
\left\{\begin{array}{ll}
1 & \mathrm{if}\ f(v)\leq i\leq g(v),\\
0 & \mathrm{otherwise},\end{array}\right.$$
for all $v\in V$ and $i\in\NN$.  In this case, if $g(v)\leq f(v)+1$
then $\Gamma_{v}(y)=\Gamma(D(v),\u^{(v)},y)$ is $\S[\pi]$--nonvanishing.
If $f(v)+1=k=g(v)-1$ then $1\leq k\leq D(v)-1$ and 
$\Gamma_{v}(y)$ has the form in Lemma 21 with $\beta=1$.
A short calculation shows that $\Gamma_{v}(y)$ is 
$\S[2\pi/3]$--nonvanishing.  Corollary 19 thus implies part (a).  For part (b)
one checks that $\Gamma_v(y)$ is $\S[3\pi/4]$--nonvanishing
if and only if either $k=1$ or $k=D(v)-1$ or $D(v)\leq 5$,
and thus Corollary 19 implies part (b).
\end{proof}

Our last application of Corollary 19 is a set of fugacities $\{\u^{(v)}\}$ weighting the
$(f,g)$--factors of a graph that is sufficient to imply logarithmic 
concavity of the resulting numbers.  This can be regarded as a (very)
weak form of Conjecture 1.
\begin{THM}
Let $G=(V,E)$ be a graph, $D=\deg(G)$, and let $f,g:V\goesto\NN$ with 
$f\leq g\leq D$.  For each $v\in V$ let $g(v)-f(v)=2a_{v}+b_{v}$ with
$0\leq b_{v}\leq  1$, and define
$$q_{v}(y)=y^{f(v)}(1+y)^{b_{v}}(1+\sqrt{3}y+y^{2})^{a_{v}}.$$
Assign nonnegative fugacities $\u^{(v)}$ to $v\in V$ so that
$\Gamma(D(v),\u^{(v)},y)$ is the result of applying Proposition
$11$\textup{(a)} to
$(1+y)^{D(v)}$ and $q_v(y)$.  Then the polynomial
$F(G;\{\u^{(v)}\},t)$ is $\S[2\pi/3]$--nonvanishing.
Therefore, the coefficients $\{N(G;\{\u^{(v)}\},j)\}$ are 
logarithmically concave.
\end{THM}
\begin{proof}
By construction, each $q_{v}(y)$ is $\S[5\pi/6]$--nonvanishing.
By Proposition 11(a), the same is true for each $\Gamma(D(v),\u^{(v)},y)$.
The result follows immediately from Corollary 19 and Proposition 14.
\end{proof}
Replacing the quadratic $1+\sqrt{3}y+y^2$ in Theorem 26 by $1+\sqrt{2}y+y^2$
gives a sufficient condition for $F(G;\{\u^{(v)}\},t)$
to be $\H$--nonvanishing.  Using $1+2y+y^2$ instead gives a sufficient
condition for $F(G;\{\u^{(v)}\},t)$ to be
$\S[\pi]$--nonvanishing.\\

Parts (b) and (c) of Theorem 17 have the following consequence.
For a graph $G=(V,E)$ and functions $f,g:V\goesto\NN$, let
$$\binom{g}{f}:=\prod_{v\in V}\binom{g(v)}{f(v)}.$$
\begin{THM}
Let $G=(V,E)$ be a graph.\\
\textup{(a)}\ Let $D=\deg(G)$ and assign nonnegative fugacities 
$\u^{(v)}$ to each $v\in V$ so that $\Gamma(D(v),\u^{(v)},y)$ has degree
$D(v)$ and only zeros of unit modulus.  Then every zero of
$F(G;\{\u^{(v)}\},t)$ has unit modulus.\\
\textup{(b)}\ In particular, every zero of the polynomial
$$\sum_{H\subseteq E}\binom{\deg(G)}{\deg(H)}^{-1}t^{\#H}$$
has unit modulus.
\end{THM}
\begin{proof}
For part (a), the hypotheses of Theorem 17(b,c) are satisfied with 
$\kappa=1$, so that $F(G;\{\u^{(v)}\},\z)$ is both $\D$--nonvanishing 
and $\E$--nonvanishing.  After diagonalizing all the variables
$z_v\mapsto t^{1/2}$, Lemma $7$ yields the result.

Part (b) is the special case of part (a) in which
the fugacities are $u^{(v)}_i=\binom{D(v)}{i}^{-1}$ for all
$v\in V$ and $0\leq i\leq D(v)$, so that the generating functions are
$\Gamma(D(v),\u^{(v)},y)=1+y+y^2+\cdots +y^{D(v)}=(1-y^{1+D(v)})/(1-y).$
\end{proof}

\section{Directed Graphs.}

We can rework the machinery of the previous section for directed graphs,
as follows.  Let $G=(V,E)$ be a directed graph.  For a spanning directed
subgraph $H\subseteq E$ of $G$, let outdeg$(H)$ and indeg$(H)$ denote the
vectors of outdegrees and of indegrees in $H$, respectively.  We use the notation
$vew\in E$ to denote that $e$ is an edge of $G$ directed out of $v$ and
into $w$.  Associate two sets of complex variables $\z':=\{z'_v:\ v\in V\}$
and $\z'':=\{z''_v:\ v\in V\}$ with the vertices of $G$.  Fix positive
real weights $\blam:=\{\lambda_e:\ e\in E\}$ for the edges of $G$.
The polynomial
$$F(G;\blam,\z',\z''):=\prod_{vew\in E}(1+\lambda_e z'_v z''_w)
=\sum_{H\subseteq E}\blam_H
(\z')^{\mathrm{outdeg}(H)}(\z'')^{\mathrm{indeg}(H)}$$
is a weighted generating function for all spanning directed subgraphs of $G$.
As in Proposition 16 this polynomial is $\H$--nonvanishing, and if 
$\blam\equiv\boldsymbol{1}$ it is also both $\D$-- and 
$\E$--nonvanishing.  For
any $D',D'':V\goesto\NN$ such that outdeg$(G)\leq D'$ and indeg$(G)\leq D''$,
each variable $z'_v$ may be polarized $D'(v)$ times, and each variable
$z''_v$ may be polarized $D''(v)$ times.  The resulting polynomial
$$P^{D',D''}F(G;\{z'_{vi}\},\{z''_{vj}\})$$
is still $\H$--, $\D$--, and $\E$--nonvanishing, by repeated application
of Proposition 10.
From this point onward, the method of proof of Theorem 17 can be applied
\emph{mutatis mutandis}, and no new complications arise.  Lacking a 
compelling application of the result, we leave the details to the reader.

\end{document}